\documentclass[11pt]{amsart}
\usepackage{amsfonts}
\usepackage{amssymb}
\usepackage[arrow,matrix]{xy}
\usepackage{amsmath,amssymb, bbm, amscd, amsthm,mathrsfs,hyperref}
\theoremstyle{plain}

\newtheorem{thm}{Theorem}[section]
\newtheorem{cor}[thm]{Corollary}
\newtheorem{lem}[thm]{Lemma}

\theoremstyle{definition}
\newtheorem{defn}[thm]{Definition}

\newtheorem{question}[thm]{Question}
\newtheorem{conjecture}[thm]{Conjecutre}
\def\vps{\varepsilon}

\def\bdt{\Delta}

\def\it{\textit}

\def\ot{\otimes}

\def\ra{\rightarrow}

\def\la{\leftarrow}

\def\Hom{\operatorname {Hom}}
\def\Ext{\operatorname {Ext}}
\def\Tor{\operatorname {Tor}}

\def\dim{\operatorname {dim}}

\def\Mod{\operatorname {Mod}}

\def\injdim{\operatorname{injdim}}

\def\kk{\mathbbm{k}}

\begin{document}
\title[Bijectivity of antipode]{\bf A note on the bijectivity of antipode of a Hopf algebra and its applications}
%
%\author{Xingting Wang}
%\address{Xingting WANG\newline Department of Mathematics, Temple University, Philadelphia, PA 19122, USA }
%\email{xingting@temple.edu}
%
%\author{Xiaolan YU}
%\address {Xiaolan YU\newline Department of Mathematics, Hangzhou Normal University, Hangzhou, Zhejiang 310036, China}
%
%\email{xlyu@hznu.edu.cn}
%
%\author{Yinhuo ZHANG}
%\address {Yinhuo ZHANG\newline Department WNI, University of Hasselt, Universitaire Campus, 3590 Diepeenbeek,Belgium } \email{yinhuo.zhang@uhasselt.be}

\author{Jiafeng L\"u}
\thanks{The first author is supported by the National Natural Science Foundation of China: No. 11571316, No. 11001245 and the Natural Science Foundation of Zhejiang Province: No. LY16A010003.}
\address{Department of Mathematics, Zhejiang Normal University, Jinhua, Zhejiang 321004, China} \email{jiafenglv@zjnu.edu.cn}
\author{Sei-Qwon Oh}
\thanks{The second author is supported by Chungnam National University Grant.}
\address{Department of Mathematics, Chungnam National  University, 99 Daehak-ro,   Yuseong-gu, Daejeon 34134, Korea} \email{sqoh@cnu.ac.kr}
\author{Xingting Wang}
\address{Department of Mathematics, Temple University, Philadelphia, 19122, USA} \email{xingting@temple.edu}
\thanks{The third author is supported by AMS-Simons travel grant.}
\author{Xiaolan Yu}
\address{Department of Mathematics, Hangzhou Normal University, Hangzhou, Zhejiang 310036, China}\email{xlyu@hznu.edu.cn}
\thanks{The fourth author is supported by the National Natural Science Foundation of China: No. 11301126, No. 11571316, No. 11671351.}

%\thanks{}
\date{}
\begin{abstract}
Certain sufficient homological and ring-theoretical conditions are given for a Hopf algebra to have bijective antipode with applications to noetherian Hopf algebras regarding their homological behaviors.
\end{abstract}

\maketitle

\section*{Introduction}\label{0}
A classical result due to Larson and Sweedler \cite{ls} states that any finite-dimensional Hopf algebra has bijective antipode.  In general, the antipode of an infinite-dimensional Hopf algebra does not need to be bijective. For instance, Takeuchi \cite{take} constructed the free Hopf algebra generated by a coalgebra whose antipode is injective but not surjective. On the other hand, Schauenburg \cite{sch} gave examples of Hopf algebras whose antipode is surjective but not injective.

In recent development, the study of infinite-dimensional Hopf algebras seems to be of growing importance, which reveals that some well-known results about finite-dimensional Hopf algebras surprisingly have incarnations in the realm of noetherian Hopf algebras (see, e.g., survey papers \cite{Bro, good}). Among these progress, it is worthy to point out that the bijectivity of the antipode frequently plays an essential role in establishing these properties (see, e.g.,\cite{bz,hoz, rrz, wyz}). Therefore, one prompts to ask for criterions concerning the bijectivity of the antipode of a Hopf algebra.

In \cite{skr}, Skryabin gave two sufficient conditions for the bijectivity, which are purely ring-theoretic. As a corollary, he proved that the antipode of any noetherian Hopf algebra is always injective, and it is surjective if certain quotient ring exits \cite[Corollary 1]{skr}. Moreover, he proposed
\begin{conjecture}(Skryabin)\label{bna}
Every noetherian Hopf algebra has bijective antipode.
\end{conjecture}

Recently, Meur showed that, by imposing a purely homological restriction, any twisted Calabi-Yau Hopf algebra has bijective antipode \cite[Proposition 1]{Meur}. The next result proved in the present paper uses both homological and ring-theoretic restrictions on a Hopf algebra.

\begin{thm}\label{Bijection}
Let $H$ be a Hopf algebra such that the left or right trivial module ${}_\vps \kk$ or $\kk_\vps$ has a resolution by finitely generated projective modules. Suppose $H$ satisfies one of the following conditions.
\begin{itemize}
\item[(i)] $\dim \Ext_{H}^i({}_\vps \kk, H)=1$ for some integer $i\ge 0$;
\item[(ii)] $\dim \Ext_{H^{\mathrm{op}}}^j(\kk_\vps,H)=1$ for some integer $j\ge 0$.
\end{itemize}
Then $H$ has injective antipode. Moreover, if both (i) and (ii) hold for $H$ with $i=j$ and $H$ additionally has one of following properties
\begin{itemize}
\item[(iii)] every left invertible element is regular;
\item[(iv)] every right invertible element is regular.
\end{itemize}
Then $H$ has bijective antipode.
\end{thm}

The class of Hopf algebras satisfies the above assumptions is large. For instance, the homological restrictions (i) and (ii) are weaker versions of AS-Gorenstein condition (see, e.g., Definition \ref{defn as}), and the ring-theoretic restrictions (iii) and (iv) are held by any Hopf algebra that is weakly finite,  which includes all noetherian Hopf algebras and Hopf domains. We are able therefore to obtain

\begin{cor}\label{App}
Any noetherian AS-Gorenstein Hopf algebra has bijective antipode.
\end{cor}

By a celebrated result of Wu and Zhang \cite{wz}, any noetherian affine PI Hopf algebra is AS-Gorenstein, which yields another proof of the following

\begin{cor}\cite[Corollary 2]{skr}\label{PI}
Any noetherian affine PI Hopf algebra has bijective antipode.
\end{cor}

Now it becomes clear that an affirmative answer to the following question \cite[Question E]{Bro} regarding the homological behaviors of noetherian Hopf algebras will help to answer Conjecture \ref{bna}.

\begin{question}(Brown)\label{Brown}
Is every noetherian Hopf algebra AS-Gorenstein?
\end{question}

The proof of our main theorem is based on analyzing the bimodule structures arising from the Hochschild cohomlogy of $H$ with coefficients in certain bimodule over $H$ (see Theorem \ref{FLemma}). With the help of Corollary \ref{App}, we apply the same idea to noetherian Hopf algebras.  We are able to extend Radford $S^4$ formula to any noetherian AS-Gorenstein Hopf algebra (see Theorem \ref{RS4}) and establish equivalent conditions regarding the homological behaviors of noetherian Hopf algebras (see Theorem \ref{Gorenstein} and Theorem \ref{CY}).

\section{Preliminaries}\label{1}
Throughout this paper, we work over a fixed field $\kk$. Unless stated otherwise all algebras and vector spaces are over $\kk$. The unadorned  tensor $\ot$ means $\ot_\kk$.
Given an algebra $A$, we write $A^\mathrm{op}$ for the \it{opposite algebra} of $A$ and $A^e$ for the \it{enveloping algebra} $A\ot A^{\mathrm{op}}$. The category of left (resp. right) $A$-modules is denoted by $\Mod(A)$ (resp. $\Mod(A^{\mathrm{op}})$). An $A$-bimodule $M$ can be identified with a left $A^e$-module, that is, an object $M$ in $\Mod(A^e)$ with action
$$(a\ot b)\cdot m=amb$$
for all $a\ot b\in A^e$ and $m\in M$. 

Note that  an $A$-bimodule $M$ can also be a right $A^e$-module with right $A^e$-action
$$m\cdot(a\ot b)=b m a$$ 
for all $a\ot b\in A^e$ and $m\in M$. Conversely, if $M$ is a right $A^e$-module then $M$ becomes an $H$-bimodule with bimodule action
 $$b m a=m\cdot(a\ot b)$$
for all $a,b\in A$ and $m\in M$.

For an $A$-bimodule $M$ and two algebra homomorphisms $\mu$ and $\nu$, we let $^\mu M^\nu$ denote the \it{twisted $A$-bimodule} such that $^\mu M^\nu\cong M$ as vector spaces, and the bimodule structure is given by
$$a\cdot m \cdot b=\mu(a)m\nu(b),$$
for all $a,b\in A$ and $m\in M$. If one of the homomorphisms is the identity, we will omit it.

We preserve $H$ for a Hopf algebra, and as usual, we use the symbols $\bdt$, $\vps$ and $S$ respectively for its comultiplication, counit, and antipode. We use  Sweedler's (sumless) notation for the comultiplication of $H$. We write ${}_\vps\kk$ (resp. $\kk_\vps$) for the left (resp. right) trivial module defined by the counit of $H$.

\begin{defn}
Let $\xi:H\ra \kk$ be an algebra homomorphism. The \textit{left winding automorphism} $\Xi^\ell_{\xi}$ of $H$ given by $\xi$ is defined to be
$$\Xi_\xi^\ell (a)=\xi(a_1)a_2,$$
for any $a\in H$. Similarly, the \textit{right winding automorphism} of $H$ given by $\xi$ is defined to be $$\Xi_{\xi}^r(a)=a_1\xi(a_2),$$
for any $a\in H$.
\end{defn}

We recall some well-known properties of winding automorphisms.

\begin{lem}(cf. \cite[Lemma 2.5]{bz})\label{windingaut}
\
\begin{enumerate}
\item[(i)] $(\Xi_\xi^\ell)^{-1}=\Xi_{\xi S}^\ell$.
\item[(ii)] $\xi S^2=\xi$, so $\Xi_{\xi}^\ell=\Xi_{\xi S^2}^\ell$.
\item[(iii)] $\Xi_\xi^\ell S^2=S^2 \Xi_\xi^\ell$.
\item[(iv)] The above are true for right winding automorphisms.
\item[(v)] Left and right winding automorphisms always commute with each other.
\end{enumerate}
\end{lem}

\begin{defn}(cf. \cite[definition 1.2]{bz})\label{defn as} Let $H$ be a noetherian Hopf algebra.
\begin{enumerate}
\item[(i)] We say $H$ has \textit{finite injective dimension} if the injective dimensions of $\!_HH$ and $H_H$ are both finite. In this case these integers are equal by \cite{za}, and
we write $d$ for the common value. We say $H$ is \textit{regular} if it has finite global dimension. Right global dimension always equals left global dimension for Hopf algebras \cite[Proposition 2.1.4]{wyz}; and, when finite, the global dimension equals the injective dimension.
\item[(ii)] The Hopf algebra $H$ is said to be \it{Artin-Schelter Gorenstein}, which we usually abbreviate to AS-Gorenstein, if
\begin{enumerate}
\item[(AS1)] $\injdim  {_HH}=d<\infty$,
\item[(AS2)] $\Ext_H^i({_\vps\kk},{H})=0$ for $i\neq d$ and $\dim\Ext_H^d({_\vps\kk},{H})=1$,
\item[(AS3)] the right $H$-module versions of (AS1,AS2) hold.
\end{enumerate}
\item[(iii)] If, in addition, the global dimension of $H$ is finite, then $H$ is called \it{Artin-Schelter regular}, which is usually shorten to AS-regular.
\end{enumerate}
\end{defn}

Suppose $H$ is noetherian AS-Gorenstein of finite injective dimension $d$. Then $\Ext^d_H({}_\vps\kk, H)$ is a one-dimensional right $H$-module. Any nonzero element in $\Ext^d_H({}_\vps\kk, H)$ is called a \it{left homological integral} of $H$. Usually, $\Ext^d_H({}_\vps\kk, H)$ is denoted by $\int^\ell_H$. Similarly, any nonzero element in $\Ext^d_{H^{op}}(\kk_\vps, H)$ is called a \it{right homological integral}. And $\Ext^d_{A^{op}}(\kk_\vps, A )$ is denoted by $\int^r_H$. Abusing the language slightly, $\int^\ell$ (resp. $\int^r$) is also called the left (resp. right) (homological) integral. Since the right $H$-module structure on $\int^\ell$ is given by some algebra homomorphism from $H$ to $\kk$, we can define left and right winding automorphisms given by $\int^\ell$. This also applies to $\int^r$ by using its left $H$-module structure. We say $H$ is \it{unimodular} if $\int^\ell\cong \kk_\vps$ as right $H$-modules. Clearly it is equivalent to the left or right winding automorphism given by $\int^\ell$ is identity.

In \cite{g2}, Ginzburg introduced Calabi-Yau algebras whose algebraic structures arises naturally in the geometry of Calabi-Yau manifolds and mirror symmetry. Calabi-Yau algebras are one of the examples satisfying the Van den Bergh duality, which was introduced by Van den Bergh \cite{vanh} in order to study Poincar\'e duality between Hochschild homology and cohomology. We adopt all these definitions to noetherian Hopf algebras.

\begin{defn}(cf. \cite{g2, vanh, bz})
Let $H$ be a noetherian Hopf algebra.
\begin{itemize}
\item[(i)] We say $H$ satisfies the \it{Van den Bergh condition} if $H$ has finite injective dimension $d$ and
\[
\Ext_{H^e}^i(H,H^e)=
\begin{cases}
0  &  i\neq d\\
U  &  i=d
\end{cases}
\]
where $U$ is an invertible $H$-bimodule. We usually call $U$ the \it{Van den Bergh dualising module} for $H$.
\item[(ii)] We say $H$ has the \it{Van den Bergh duality} if it satisfies the Van den Bergh condition and $H$ is \it{homologically smooth}, that is, $H$ has a bounded resolution in $\Mod(H^e)$ by finitely generated projective modules.
\item[(iii)] We say $H$ is \it{twisted Calabi-Yau} if $H$ has the Van den Bergh duality with the Van den Bergh dualising module given by $H^\nu$ for some algebra automorphism $\nu$ of $H$. Moreover, we say $H$ is \it{Calabi-Yau} if $\nu$ can be chosen as an inner automorphism.
\end{itemize}
\end{defn}

\section{An isomorphism lemma for Hopf bimodules}
In this section, we aim at investigating the bimodule structures arising from the Hochschild cohomology of $H$ with coefficients in the envelop algebra $H^e$. In particular, we do not require $H$ to be noetherian or have bijective antipode.

Note that the following map $$(1\ot S)\Delta: H\to H^e,\ \ a\mapsto a_1\ot S(a_2)$$ is an algebra homomorphism.

\begin{defn}
 We define the left adjoint functor $\mathscr L$ from the category of left $H^e$-modules into the category of left $H$-modules such that, for every left $H^e$-module $M$,  $\mathscr L(M)=M$ as vector spaces with the left action
$$a\cdot m=(1\ot S)\Delta(a)\cdot m=(a_1\ot S(a_2))\cdot m $$
for $a\in H$ and $m\in M$.
Similarly, the right adjoint functor $\mathscr R$ from the category of right $H^e$-modules into the category of right $H$-modules such that, for every right $H^e$-module $M$, $\mathscr R(M)=M$ as vector spaces with the right action
$$ m\cdot a= m\cdot (1\ot S)\Delta(a)=m\cdot(a_1\ot S(a_2))$$
for $a\in H$ and $m\in M$.
\end{defn}

Here we introduce natural module actions and elementary properties which will be used.
Since the envelope algebra $H^e$ is an algebra, $H^e$ is equipped with a natural $H^e$-bimodule structure induced by the multiplication of $H^e$. That is, the  left action  is given by
\begin{equation}\label{HHlm}
(a\ot b)\ra (x\ot y)=(a\ot b)(x\ot y)=a  x\ot y  b,\end{equation}
called the \it{outer action},
and the right action is given by
\begin{equation}\label{HHrm}(x\ot y)\la(a\ot b)=(x\ot y)(a\ot b) =x  a\ot b  y,\end{equation}
called the \it{inner action}.
As a consequence, $\mathscr L(H^e)$ can be viewed as an $H$-$H^e$-bimodule, where the left $H$-action is given by applying the left adjoint functor to the outer action
$$a\cdot (x\ot y)=((1\ot S)\Delta(a))(x\ot y)=a_1  x\ot y  S(a_2)$$
and the inner action gives the right $H^e$-module structure. On the other hand, $\mathscr R(H^e)$ is an $H^e$-$H$-bimodule with the right action
$$(x\ot y)\cdot a=(x\ot y)((1\ot S)\Delta(a))=xa_1\ot S(a_2)y$$
together with the outer action for the left $H^e$-module structure.

Let $M$ and $N$ be two left $H$-modules. Then $M\ot N$ is a left $H\ot H$-module with a natural left $H\ot H$-action
$$(a\ot b)\ra (x\ot y)=(a\cdot  x)\ot (b\cdot y).$$   
Since there are two natural algebra homomorphisms from $H$ into $H\ot H$ such that 
$$H\to H\ot H,\ \ \ a\mapsto a\ot 1$$
and 
$$H\to H\ot H,\ \ \ a\mapsto 1\ot a,$$
there are two left $H$-module actions on $M\ot N$ such that
$$a\cdot(x\ot y)=(a\ot 1)\ra (x\ot y)=(a\cdot x)\ot y,\ \ \ (\text{denoted by ${}_*M\ot N$})$$
and
$$a\cdot(x\ot y)=(1\ot a)\ra (x\ot y)=x\ot (a\cdot y). \ \ \ (\text{denoted by $M\ot {}_*N$})$$
Analogously, for any right $H$-modules $M$ and $N$, there are two right $H$-module actions $M_*\ot N$ and $M\ot N_*$.

Since the co-multiplication map $\Delta:H\to H\ot H$ is an algebra homomorphism, every left (respectively, right) $H\otimes H$-module becomes a left
(respectively, right) $H$-module with the action induced by $\Delta$, namely
$$a\cdot(x\ot y)=\Delta(a)\ra (x\ot y)=(a_1\cdot x)\ot (a_2\cdot y).$$

Let $R$ and $T$ be algebras. For a left $R$-module $_RN$ and an $R$-$T$-bimodule $_RM_T$, $\text{Hom}_R(_RN,_RM_T)$ is a right $T$-module
with the right $T$-action
$$(ft)(n)=f(n)t$$
for $f\in \text{Hom}_R(_RN,_RM_T)$, $t\in T, n\in N$. For a right $T$-module $N_T$ and an $R$-$T$-bimodule $_RM_T$, $\text{Hom}_T(N_T,_RM_T)$ is a left $R$-module
with the left $R$-action
$$(rf)(n)=rf(n)$$
for $f\in \text{Hom}_T(N_T,_RM_T)$, $r\in R, n\in N$.
We often write $\text{Hom}_{T^{\text{op}}}(N_T,_RM_T)$ for $\text{Hom}_T(N_T,_RM_T)$.
For a $R$-$T$-bimodule $_RN_T$ and a left  $R$-module $_RM$, $\text{Hom}_R(_RN_T,_RM)$ is a left $T$-module
with the left $T$-action
$$(tf)(n)=f(nt)$$
for $f\in \text{Hom}_R(_RN_T,_RM)$, $t\in T, n\in N$.

The following is parallel to Lemma 2.4 in \cite{bz} and Lemma 2.1.2 in \cite{wyz}. For the sake of completeness, we include a proof here.

\begin{lem}\label{adjointH} Let $A$ be an algebra. There are natural isomorphisms for all integers $i\ge 0$.
\begin{enumerate}
\item[(i)] Let $M$ be an $H^e$-$A$-bimodule. Then $\Ext^i_{H^e}(H,M)\cong \Ext^i_H({}_\vps\kk,\mathscr L(M))$ as right $A$-modules.
\item[(ii)] Let $M$ be an $A$-$H^e$-bimodule. Then  $\Ext^i_{H^e}(H,M)\cong \Ext^i_{H^\mathrm{op}}(\kk_\vps,\mathscr R(M))$ as left $A$-modules.
\end{enumerate}
\end{lem}
\proof We only prove (i), the proof of (ii) is quite similar. Note that $H^e$-$A$-bimodule $N$ is canonically a left $H^e\otimes A^\text{op}$-module and that $H^e\ot A^\text{op}$ is a right $H^e$-module with the right action induced by the
multiplication of  $H^e\ot A^\text{op}$ since $H^e$ is considered as a subalgebra of $H^e\ot A^\text{op}$ by the inclusion map
$H^e\to H^e\ot A^\text{op}$, $x\mapsto x\ot1$.
 First of all, one sees easily that any injective $H^e$-$A$-bimodule $N$ is still injective when viewed as a left $H^e$-module since
\begin{align*}
\Hom_{H^e}(-,N)&\, \cong \Hom_{H^e}(-, \Hom_{H^e\otimes A^{\mathrm{op}}}((H^e\otimes A^{\mathrm{op}})_{H^e},N))\\
&\, \cong \Hom_{H^e\otimes A^\text{op}}((H^e\otimes A^{\mathrm{op}})_{H^e}\otimes -,N)
\end{align*}
by \cite[Theorem 2.11]{Rot}.

%The $H^e$-$A$-bimodule N can be viewed as an He ? Bop-module. It can be embedded into an injective He ? Bop-module I. We have
%HomHe (?, I) ?= HomHe (?, HomHe?Bop ((He ? Bop)He , I)) ?= HomHe?Bop ((He ? Bop)He ? ?, I).
%He ?Bop is clearly free as a He-module. Therefore, the functor HomHe (?, I) is exact. That is, I is injective as an He-module.

Next, we view $H^e$ as an $H^e$-$H$-bimodule, where the left $H^e$-action is given by \eqref{HHlm} and the right $H$-action is given by
$$(x\ot y)\cdot a=xa_1\ot S(a_2)y.$$
We simply denote it as ${}_{H^e}H^e_H$, which is free as a right module by the fundamental theorem of Hopf modules. Indeed,
there is an
$H^e$-$H$-bimodule isomorphism ${}_{H^e}H^e_H \ra H_*\ot H$  defined by $\;x\ot y \mapsto x_1 \ot x_2 y$ with inverse given by $x\ot y \mapsto x_1 \ot S(x_2) y$, where the left $H^e$-action on $H_*\ot H$ is given by $(a\ot b)\cdot(x\otimes y)=a_1x\ot a_2yb$ and the right $H$-action on $H_*\ot H$ is given by $(x\ot y)\cdot a=(x\ot y)(a\ot 1)=xa\ot y$. Since $\mathscr L\cong \Hom_{H^e}({}_{H^e}H^e_H,-)$ as functors, one gets that
\begin{align*}
\Hom_H(-,\mathscr L(M))&\cong \Hom_H(-,\Hom_{H^e}({}_{H^e}H^e_H,M))\\
&\cong \Hom_{H^e}({}_{H^e}H^e_H\ot_H-,M).
\end{align*}
As a consequence, $\mathscr L$ is exact and preserves injectivity.

Since there is an isomorphism ${}_{H^e}H^e_H \ra H_*\ot H$ by the above paragraph, we have the canonical isomorphism
$${}_{H^e}H^e_H\ot_H {}_\vps\kk\cong H_*\ot H\ot_H {}_\vps\kk\cong {}_{H^e}H.$$
 Hence (i) holds for $i=0$.
It follows that (i) holds for all $i\ge 0$ by taking an injective resolution of $M$ as $H^e$-$A$-bimodules.\qed

%It is important to point out that in the above lemma if $M$ is equipped with an $H^e$-$H^e$-bimodule structure, then the isomorphisms of Ext-groups are respect to the remaining $H^e$-module structures.

\begin{lem}\label{HomProj}
\
\begin{itemize}
\item[(i)] Let $P$ be a finitely generated projective left $H$-module. Then
$$\Hom_H(P,\mathscr L (H^e))\cong \Hom_H(P,H)\otimes \!_*H^{S^2}$$
as $H$-bimodules, where the bimodule structure on $\Hom_H(P,H)\otimes \!_*H^{S^2}$ is given by $a(x\otimes y)b=xb_1\otimes ayS^2(b_2)$.
\item[(ii)] Let $Q$ be a finitely generated projective right $H$-module. Then
$$\Hom_{H^\mathrm{op}}(Q,\mathscr R(H^e))\cong \!^{S^2}H_*\otimes \Hom_H(P,H)$$
as $H$-bimodules, where the bimodule structure on $\!^{S^2}H_*\otimes \Hom_H(P,H)$ is given by $a(x\otimes y)b=S^2(a_1)xb\otimes a_2y$.
\end{itemize}
\end{lem}
\proof
(i) Note that $\Hom_H(P,\mathscr L (H^e))$ is a left $H^e$-module and thus a $H$-bimodule since $\mathscr L (H^e)$ is a $H$-$H^e$-bimodule. First of all, we claim that $\mathscr L(H^e)\cong \!_*H\otimes H:=V$ as $H$-$H^e$-bimodules, where the left $H$-action on $V$ is defined by the left multiplication on the first factor $H$ of $V$ and the right $H^e$-action is given by $(x\otimes y)\la(a\otimes b)=xa_1\otimes byS^2(a_2)$. It can be proved via the explicit $H$-$H^e$-isomorphism $\mathscr L(H^e)\to V$ defined by $x\otimes y\mapsto x_1\otimes yS^2(x_2)$ with inverse given by $x\otimes y\mapsto x_1\otimes yS(x_2)$.

Next for any left $H$-module $M$, there exists a natural $H$-bimodule map
$$\Phi_M: \Hom_H(M,H)\otimes \!_*H^{S^2}\to \Hom_H(M,\mathscr L(H^e))\cong \Hom_H(M,V)$$
defined by $\Phi_M(f\otimes h)(m)=f(m)\otimes h$. One checks that $\Phi$ commutes with finite direct sum, that is, $\Phi_{\oplus_{i\in I} M_i}=\bigoplus_{i\in I} \Phi_{M_i}$ since the following diagram
\[
\xymatrix{
\Hom_H(\bigoplus_{i\in I} M_i, H)\otimes \!_*H^{S^2} \ar[rr]^-{\Phi_{\oplus_{i\in I}M_i}}\ar[d]^-{\cong} &&  \Hom_H(\bigoplus_{i\in I}M_i, V)\ar[d]^-{\cong}\\
\bigoplus_{i\in I} \left(\Hom_H(M_i,H)\otimes \!_*H^{S^2}\right) \ar[rr]^-{\bigoplus_{i\in I} \Phi_{M_i}} && \bigoplus_{i\in I} \Hom_H(M_i,V)
}
\]
commutes whenever $I$ is a finite index set. Suppose $P$ is finitely generated projective, then there exists another left $H$-module $Q$ such that $P\bigoplus Q=\bigoplus_{i\in I}H_i$ over a finite index set $I$, where each $H_i\cong H$ as left $H$-modules. Note that $\Phi_H$ is clearly an isomorphism. Hence $\Phi_P\bigoplus \Phi_Q=\Phi_{P\oplus Q}=\bigoplus_{i\in I} \Phi_{H_i}$ is an isomorphism, which implies that $\Phi_P$ is an isomorphism.

Finally, denote by $W=H\otimes H_*$ the $H^e$-$H$-bimodule, where the right $H$-action is the right multiplication on the second factor $H$ of $W$ and the left $H^e$-action is given by $(a\otimes b)\ra (x\otimes y)=S^2(a_1)xb\otimes a_2y$. Then (ii) can be proved in the same fashion by using the $H^e$-$H$-isomorphism $\mathscr R(H^e)\cong W$ via $x\otimes y\mapsto S^2(x_1)y\otimes x_2$ with inverse $x\otimes y\mapsto y_2\otimes S(y_1)x$.
\qed

\begin{lem}\label{Compact}
The following are equivalent.
\begin{itemize}
\item[(i)] $H$ has a resolution in $\Mod(H^e)$ by finitely generated projective modules.
\item[(ii)] ${}_\vps\kk$ has a resolution in $\Mod(H)$ by finitely generated projective modules.
\item[(iii)] The right $H$-module version of (ii) holds.
\end{itemize}
\end{lem}
\proof
(i) $\Rightarrow$ (ii), (iii) 
Let $M$ be an $H$-bimodule and let $I=\ker\vps$. Then it is easy to see that $\kk_\vps\otimes_H M\cong M/IM$.
Let $\mathcal B^\bullet$ be a resolution of $H$ in $\Mod(H^e)$ by finitely generated projective modules. Then, using the above result, one can observe easily that $\kk \otimes_H \mathcal B^\bullet$ is a resolution of ${}_\vps\kk_\vps  \otimes_H H\cong {}_\vps\kk$ in $\Mod(H)$ by finitely generated projective modules. It is the same for (iii) when we tensor $\otimes_H{}_\vps\kk_\vps$ on the right side of $\mathcal B^\bullet$.

(iii), (ii) $\Rightarrow$ (i)  In the proof of Lemma \ref{adjointH}, one sees that the left adjoint functor $\mathscr L: \Mod(H^e)\to \Mod(H)$ is just a restriction functor, which certainly commutes with direct limits. Applying \cite[Corollary P. 130]{bro}, $\Ext_H^i({}_\vps\kk,\mathscr L(-))$ commutes with direct limits for all $i$ for ${}_\vps\kk$ has a resolution in $\Mod(H)$ by finitely generated projective modules. This implies that $\Ext_{H^e}^i(H,-)$ commutes with direct limits in $\Mod(H^e)$ for all $i$ since $\Ext_{H^e}^i(H,-)\cong \Ext_H^i({}_\vps\kk,\mathscr L(-))$ by Lemma \ref{adjointH}. Then one concludes again by \cite[Corollary P. 130]{bro} that $H$ has a resolution in $\Mod(H^e)$ by finitely generated projective modules. The proof for (iii) is exactly the same.

\qed

\begin{thm}\label{FLemma}
Assume the conditions in Lemma \ref{Compact} hold. Then there are $H$-bimodule isomorphisms
$$\Ext_{H^e}^i(H,H^e)\cong \Ext_H^i\left({}_\vps\kk,H\right)\otimes \!_*H^{S^2}\cong\, \!{}^{S^2}H_* \otimes \Ext^i_{H^\mathrm{op}}\left(\kk_\vps,H\right)$$
for all $i$, where the bimodule structures on the second and third ones are given by $a(x\otimes y)b=xb_1\otimes ayS^2(b_2)$ and $a(x\otimes y)b=S^2(a_1)xb\otimes a_2y$, respectively.
\end{thm}
\proof
Since $(H^e)^\mathrm{op}\cong H^e$, there is an equivalence between the category of left $H^e$-modules and the category of right $H^e$-modules. As a consequence, $\Ext_{H^e}^i(H,H^e)$ can be computed  by using both the outer action and the inner action of $H^e$ defined in \eqref{HHlm} and \eqref{HHrm}, respectively.

First of all, we use the outer action \eqref{HHlm} on $H^e$ to compute the Hochschild cohomology $\Ext_{H^e}^i(H,H^e)$. By Lemma \ref{Compact}, we can take $\mathcal P^\bullet$ to be a resolution of ${}_\vps\kk$ in $\Mod(H)$ consisting of finitely generated projective modules. Then we have
\begin{align*}
\Ext_{H^e}^i(H,H^e)&\, =\Ext_{H}^i\left({}_\vps \kk, \mathscr L(H^e)\right) \tag{Lemma \ref{adjointH}}\\
&\,=\mathrm{H}^i(\Hom_H(\mathcal P^\bullet, \mathscr L(H^e))\\
&\,=\mathrm{H}^i(\Hom_H(\mathcal P^\bullet, H)\otimes \!_* H^{S^2})  \tag{Lemma \ref{HomProj}}\\
&\,=\mathrm{H}^i(\Hom_H(\mathcal P^\bullet, H))\otimes \!_*H^{S^2} \tag{K\"unneth formula}\\
&\,=\Ext_H^i\left({}_\vps\kk, H\right)\otimes \!_*H^{S^2}.
\end{align*}

On the other hand, we can apply the inner action \eqref{HHrm} on $H^e$ to compute the Hochschild cohomology $\Ext_{H^e}^i(H,H^e)$. We get $\Ext_{H^e}^i(H,H^e)\cong\!^{S^2}H_* \otimes \Ext^i_{H^\mathrm{op}}(\kk_\vps,H)$ by the same argument. This proves the result.
\qed

\proof[Proof of Theorem \ref{Bijection}]
For the injectivity of $S$, suppose (i) holds for $H$ and the proof for (ii) is analogous. Note that $\text{Hom}_{H}(M,H)$
is a right $H$-module for any left $H$-module $M$. Hence we can write $\Ext_{H}^i({}_\vps\kk,H)=\kk^\xi$ for some $\xi\in \Hom_{\mathrm{Alg}}(H,\kk)$. For simplicity, we denote the left winding automorphism $\Xi_\xi^\ell$ still by $\xi$. By Theorem \ref{FLemma}, we have the following isomorphisms
\begin{align}\label{BiExt}
\!^{S^2}H_*\otimes \Ext_{H^\mathrm{op}}^i(\kk_\vps,H)\cong \Ext_H^i({}_\vps\kk, H)\otimes \!_* H^{S^2}\cong \kk^\xi\otimes\!_*H^{S^2}\cong H^{S^2\xi}
\end{align}
as $H$-bimodules. Since the very left side of \eqref{BiExt} is a free right $H$-module, this implies that $H^{S^2\xi}$ is torsion free on the right side. Thus $S$ is injective.

Now assume that (i) and (ii) both hold for $H$ with $i=j$. Then we can further write $\Ext_{H^\mathrm{op}}^i(\kk_\vps,H)=\!^\eta\kk$ for some $\eta\in \Hom_{\mathrm{Alg}}(H,\kk)$. We still denote by $\eta$ the right winding automorphism $\Xi_\eta^r$. Then it is straightforward to check that \eqref{BiExt} implies that $\!^{S^2\eta}H$ and $H^{S^2\xi}$ are isomorphic as $H$-bimodules. Take $\Phi: \!^{S^2\eta}H\to H^{S^2\xi}$ to be such an isomorphism with inverse $\Phi^{-1}$. Denote by $x=\Phi(1)$ and $y=\Phi^{-1}(1)$. One immediately, by the definition of the inverse $\Phi\Phi^{-1}=\mathrm{id}=\Phi^{-1}\Phi$, verifies that the following hold in $H$ for any $a,b\in H$.
\begin{align}\label{S4}
xS^4\xi\eta(a)S^2\xi(y)=a,\quad S^2\eta(x)S^4\xi\eta(b)y=b.
\end{align}
Here we use the fact that $\xi,\eta, S^2$ commute with each other by Lemma \ref{windingaut}. Let $a=b=1$. One gets $xS^2\xi(y)=S^2\eta(x)y=1$.

If (iii) holds. Applying $S^2\eta$ to $xS^2\xi(y)=1$, one sees that $S^2\eta(x)S^4\xi\eta(y)=S^2\eta(x)y=1$. So $S^2\eta(x)(y-S^4\xi\eta(y))=0$, which implies that $y=S^4\xi\eta(y)$ since $S^2\eta(x)$ is left invertible hence it is not a left zero divisor by (iii). Then one gets
\begin{align*}
S^2\xi(y)x=S^2\xi(y)\ra\Phi(1)=\Phi(S^2\xi(y)\ra 1)=\Phi(S^4\xi\eta(y))=\Phi(y)\\
=\Phi(1\la y)=\Phi(1)\la y=x S^2\xi(y)=1.
\end{align*}
Thus $x$ and $S^2\xi(y)$ are invertible to each other. One obtains from \eqref{S4} the following formula
\begin{align}\label{ProofS4}
S^4\xi\eta(a)=S^2\xi(y)ax.
\end{align}
As a consequence, $S^4\xi\eta$ is an inner automorphism given by the conjugation of the element $x$. Thus $S$ is bijective. Finally, the argument for (iv) is similar. This proves the result.
\qed

\section{Applications to noetherian Hopf algebras}
In this section, we apply our result to noetherian Hopf algebras satisfying the AS-Gorenstein condition, which now we know have bijective antipodes by Corollary \ref{App}. We refine many results focusing on their homological behaviors, some of which were originally stated with the assumption of the bijectivity of the antipode (see, e.g., \cite{bz,hoz}). The first result is known to be the generalization of the famous Radford's $S^4$ formula \cite{rad} to the noetherian AS-Gorenstein Hopf algebra case by Brown and Zhang. We give another proof based on Theorem \ref{Bijection}.

\begin{thm}\cite[Corollary 4.6]{bz}\label{RS4}
Let $H$ be a noetherian AS-Gorenstein Hopf algebra. Then
$$S^4=\gamma \circ \phi \circ \xi^{-1}$$
where $\xi$ and $\phi$ are respectively the left and right winding automorphisms given by the left integral of $H$, and $\gamma$ is an inner automorphism.
\end{thm}
\proof
The result basically can be derived from the proof of Theorem \ref{Bijection}. First of all, one checks that all the assumptions in Theorem \ref{Bijection} are satisfied when $H$ is noetherian AS-Gorenstein. Namely, noetherianness guarantees that ${}_\vps\kk$ admits a resolution in $\Mod(H)$ by finitely generated projective modules. Conditions (i) and (ii) follow from AS-Gorenstein assumption with $i=j=d$. Note that in a noetherian ring, a left or right invertible element is always invertible and hence it is regular (cf. \cite[Exercise 5ZE]{gw}). So (iii) and (iv) hold.

Now we keep the same notations as in the proof of Theorem \ref{Bijection}. Denote by $\xi$ the left winding automorphism given by the left integral $\int^\ell$ and $\eta$ the right winding automorphism given by the right integral $\int^r$. We write $\int^r=\!^\pi \kk$ for some $\pi\in \Hom_{\mathrm{Alg}}(H,\kk)$. By \cite[Lemma 2.1]{lwz} (note that $S$ is bijective), $\int^l=S(\int^r)=k^{\pi S}$. So by using Lemma \ref{windingaut}, one sees that $\eta^{-1}=(\Xi_\pi^r)^{-1}=\Xi_{\pi S}^r:=\phi$ is the right automorphism given by $\int^l$.

Finally, from \eqref{ProofS4} one gets that $S^4\eta\xi$ is an inner automorphism of $H$, which we now denote by $\gamma$. Note that $S^4, \eta, \xi$ and $\gamma$ commute with each other. Hence $S^4=\gamma \circ \eta^{-1} \circ \xi=\gamma \circ \phi \circ \xi^{-1}$.
\qed

\begin{question} (Brown-Zhang)
What is the inner automorphism $\gamma$ in Theorem \ref{RS4}?
\end{question}

The answer is known when $H$ is finite-dimensional, where $\gamma$ is the conjugation by the distinguish group-like element of $H$ given by $\int^\ell_{H^*}$. In view of Question \ref{Brown}, we expect Theorem \ref{RS4} should hold for any noetherian Hopf algebra.

Next, we establish several equivalent conditions regarding noetherian AS-Gorenstein and AS-regular Hopf algebras. Recall that the noncommutative version of the daulising complex was first introduced by Yekutieli in \cite{ye}, and rigid dualising complex was later introduced by Van den Bergh in \cite{vdb} in order to remedy its uniqueness.

\begin{thm}\label{Gorenstein}
Let $H$ be a noetherian Hopf algebra. Then the following are equivalent.
\begin{itemize}
\item[(i)] $H$ is AS-Gorenstein.
\item[(ii)] $H$ satisfies Van den Bergh condition.
\item[(iii)] $H$ has a rigid dualising complex $R=V[s]$, where $V$ is invertible and $s\in \mathbb Z$.
\end{itemize}
In these cases, the rigid dualising complex is $R=\!^{S^2\xi}H[d]$, where $\xi$ is the left winding automorphism given by the left integral of $H$ and $d$ is the injective dimension of $H$.
\end{thm}
\proof
(ii) $\Leftrightarrow$ (iii) follows from \cite{vdb}, also see \cite[Proposition 4.3]{bz}.

(i) $\Rightarrow$ (iii) is \cite[Proposition 4.5]{bz}, where the assumption of the bijectivity of the antipode is automatically satisfied with the help of Corollary \ref{App}.

(ii) $\Rightarrow$ (i) Suppose $H$ satisfies Van den Bergh condition with injective dimension $d$. In view of Theorem \ref{FLemma}, one sees that $\Ext_H^i({}_\vps\kk,H)=0$ for $i\neq d$ and $\Ext_H^d({}_\vps\kk,H)\neq 0$. This holds for the right side versions of the Ext-groups as well. Moreover, the Van den Bergh dualising module $U$ is isomorphic to $\Ext_H^d({}_\vps\kk,H)\otimes \!_*H$ as $H$-bimodules, where the latter one is a free left $H$-module with basis given in $\Ext_H^d({}_\vps\kk,H)$. Since $U$ is invertible, it is finitely generated projective when viewed as a left $H$-module. It can be verified by considering the autoequivalence functor $U\otimes_A-: \Mod(H)\to \Mod(H)$ with inverse functor given by $U^{-1}\otimes _A-$. Note that a left $H$-module $M$ is finitely generated if and only if $\Hom_H(M,-)$ commutes with inductive direct limits, which is certainly preserved under any autoequivalence functor. Hence $U=U\otimes_A A$ is finitely generated. As a consequence, this implies that $\Ext_H^d(\kk,H)$ is finite-dimensional. By the same reason, $\Ext_{H^\mathrm{op}}^d(\kk_\vps,H)$ is finite-dimensional. Then \cite[Lemma 3.2]{bz} shows that $H$ is AS-Gorenstein.

Finally, the formula of the rigid dualising complex is given in \cite[Proposition 4.5]{bz}.
\qed

\begin{thm}\label{CY}
Let $H$ be a noetherian Hopf algebra. Then the following are equivalent.
\begin{itemize}
\item[(i)] $H$ is twisted Calabi-Yau.
\item[(ii)] $H$ has Van den Bergh duality.
\item[(iii)] $H$ is AS-regular.
\item[(iv)] $H$ is regular and $\Ext_H^i({}_\vps\kk,H)$ are finite-dimensional for all $i$.
\item[(v)] $H$ is regular and $\Ext_{H^\mathrm{op}}^i(\kk_\vps,H)$ are finite-dimensional for all $i$.
\end{itemize}
Moroever, $H$ is Calabi-Yau if and only if $H$ is unimodular and $S^2$ is inner.
\end{thm}
\proof
(i) $\Leftrightarrow$ (ii) follows from \cite[Theorem 3.5.1]{Meur}.

(iii) $\Rightarrow$ (iv), (v) are clear.

(i) $\Leftrightarrow$ (iii)  can be easily deduced from Theorem \ref{Gorenstein} ((i) $\Leftrightarrow$ (ii)). Since if $H$ is noetherian, then it is regular if and only if it is homologically smooth. One direction is clear. The other direction: suppose $H$ is regular, then ${}_\vps\kk$ has a bounded resolution in $\Mod(H)$ by finitely generated projective modules. This implies that $H$ is homologically smooth by \cite[Proposition 2.1.5]{wyz}.

It remains to show that (iv), (v) $\Rightarrow$ (iii). Here we only prove (iv) $\Rightarrow$ (iii) and the other one is similar. We will use Ischebeck's spectral sequence such that
\[
\Ext_{H^\mathrm{op}}^p(\Ext_H^{-q}({}_\vps\kk,H), H)\Rightarrow \Tor^H_{-p-q}=\begin{cases}  \kk & p+q=0 \\ 0  & \text{elsewhere} \end{cases}.
\]
By \cite[Proposition 2.1.4]{wyz}, the right and left global dimension of $H$ are all equal to $d$. Since $H$ is noetherian, one sees that $\Ext_H^d({}_\vps\kk,H)\neq 0$ and $\Ext_{H^\mathrm{op}}^d(\kk_\vps,H)\neq 0$ \cite[\S 1.12]{brg}. Applying \cite[Proposition 1.3]{brg}, we have
$$
\Ext_{H^\mathrm{op}}^d(\Ext_H^i({}_\vps\kk,H), H)\cong \Ext_{H^\mathrm{op}}^d(\kk_\vps, H)^{\oplus \dim \Ext_H^i({}_\vps\kk,H)},
$$
where we use the fact that $\Ext_H^i({}_\vps\kk,H)$ are finite-dimensional for all $i$. This implies that $\Ext_H^i({}_\vps\kk,H)=0$ for all $i \neq d$ since the lowest degree nonvanishing term $\Ext_H^i({}_\vps\kk,H)\neq 0$ if exists with $i<d$ would contribute to the final page of the spectral sequence at the highest degree $d-i\neq 0$. In view of Theorem \ref{FLemma}, one gets $\Ext_{H^\mathrm{op}}^i(\kk_\vps,H)=0$ for all $i \neq d$. Finally, a dimension argument used in \cite[Lemma 3.2]{bz} yields that $\dim \Ext_H^d({}_\vps\kk,H)=\dim \Ext_{H^\mathrm{op}}^d(\kk_\vps,H)=1$. This proves that $H$ is AS-Gorenstein of injective dimension $d$ and hence AS-regular.

Finally, the Calabi-Yau property is \cite[Theorem 2.3]{hoz}.
\qed

\vspace{0.2in}
\noindent
{\bf Acknowledgments}
The second, third and fourth authors are grateful for the hospitality of the first author at Zhejiang Normal University summer 2016 during the time the project was started.  And the third author wants to thank James Zhang and Martin Lorenz for some helpful conversations.

%
%
%
%
%
%\subsection*{Acknowledgement}The second author is supported by grants from NSFC (No. 11301126, No. 11571316, No. 11671351).
%
\vspace{5mm}

%%%%%%%%%%%%%%%%%%%%
%\bibliographystyle{amsplain}
%\bibliography{reference}
%%%%%%%%%%%%%%%%%%%%%%%%%%%%%%%%%%%

\bibliography{}

\end{document}